\documentclass[12pt,a4paper,leqno,twoside]{amsart}

\usepackage{latexsym}
\usepackage{amsmath}
\usepackage{amsfonts}
\usepackage{amssymb}

\numberwithin{equation}{section}

\theoremstyle{plain}
\newtheorem{thm}{Theorem}[section]
\newtheorem{cor}[thm]{Corollary}
\newtheorem{lemma}[thm]{Lemma}
\newtheorem{prop}[thm]{Proposition}

\theoremstyle{definition}

\newtheorem{remark}[thm]{Remark}

\newtheorem*{ack}{Acknowledgement}

\input xy
\xyoption{all}

\def\Zz{{\mathbb Z}}
\def\Rr{{\mathbb R}}
\def\Cc{{\mathbb C}}
\def\Hh{{\mathbb H}}
\def\Ff{{\mathbb F}}
\def\Kk{{\mathbb K}}
\def\O{O}

\def\SO{SO}
\def\U{U}
\def\SU{SU}
\def\Sp{Sp}
\def\Spin{Spin}

\begin{document}

\title{$G$-structures on spheres}

\author[M.~\v Cadek and M.~C.~Crabb]{Martin \v Cadek and Michael Crabb}


\address{\newline Department of Algebra and Geometry, Masaryk University,
Jan\' a\v ckovo n\' am. 2a, 602 00 Brno, Czech Republic}
\email{cadek@math.muni.cz}

\address{\newline Department of Mathematical Sciences, University
of Aberdeen, Aberdeen AB24 3UE, U.K.}
\email{m.crabb@maths.abdn.ac.uk}

\date {October 4th, 2005}

\thanks{Research of the first author supported by the grant
MSM 0021622409 of the Czech Ministry of Education.}



\abstract {A generalization of classical theorems on the existence of sections
of real, complex and quaternionic Stiefel manifolds over spheres
is proved. We obtain a complete list of Lie group
homomorphisms $\rho:G\to G_n$, where $G_n$ is one of the groups 
$\SO(n),\, \SU(n),\, \Sp(n)$ and 
$G$ is one of the groups $\SO(k),\, \SU(k),\, \Sp(k)$, 
which reduce the structure group $G_n$
in the fibre bundle $G_n\to G_{n+1}\to G_{n+1}/G_n$.}
\endabstract

\maketitle

\section{Introduction}

Consider the fibrations
\begin{align}\label{pf1}
\SO (n) \to  \SO (n+1) & \to \SO (n+1)/\SO (n)\, = \, S^{n} \, , \\
\SU (n) \to  \SU (n+1) & \to \SU (n+1)/\SU (n)\, = \, S^{2n+1} \, , \\
\Sp (n) \to  \Sp (n+1) &\to \Sp (n+1)/\Sp (n)\,  = \, S^{4n+3} \, .
\end{align}
To deal with the three cases we shall write
$G_n = \SO (n)$, $\SU (n)$ or $\Sp (n)$ and $d=1$, $2$ or $4$ as
appropriate.
These principal bundles reduce the structure group $\O (d(n+1)-1)$
of the tangent bundle of the sphere $G_{n+1}/G_{n}=S^{d(n+1)-1}$
to the subgroup $G_n$.

Let $G$ be a Lie group and $\rho:G\to G_n$ a homomorphism. We say that
the structure group $G_n$ of a principal fibre bundle
$\tau$ over a CW-complex $X$ can be reduced to $(G,\rho)$ if the
classifying map of this fibration $\tau:X\to BG_n$ can
be factored (up to homotopy) through $B\rho:BG\to BG_n$.
$$
\xymatrix{
& BG \ar[d]^{B\rho} \\
X \ar[r]_\tau \ar[ur] & BG_n
}
$$

This paper deals with the problem of determining those groups $G$ and
homomorphisms $\rho$ to which the structure group $G_n$ in the fibrations
above can be reduced. The problem has been solved in many interesting
special cases. Considering standard inclusions $\rho:G=G_k\hookrightarrow
G_n$ we get the famous problem on sections of Stiefel manifolds over spheres
resolved in \cite{Adams1},
\cite{AW}, \cite{AT} and \cite{SS}. The other standard inclusions
$\SU(k)\hookrightarrow \SO(n)$, $\Sp(k)\hookrightarrow \SO(n)$ and 
$\Sp(k)\hookrightarrow \SU(n)$ are dealt with in \cite{D},
\cite{O1} and \cite{O2}, respectively.
In these cases the question was to
find a minimal standard subgroup to which $G_n$ can be reduced.

In \cite{L} Leonard asked an opposite question: find all maximal
proper subgroups to which $G_n$ can be reduced. He solved it
in the cases when $G$ is a reducible maximal
subgroup of $G_n$. Moreover, he proved that $G_n$ cannot be reduced 
to any proper subgroup $(G,\rho)$ if

\begin{enumerate}
\item $n$ is even and $G_n=\SO(n)$ or $\SU(n)$, unless $G_n=\SO(6)$ and
$G=\SU(3)$;
\item $n\not\equiv 11\mod 12$ and $G_n=\Sp(n)$;
\item $G$ is a nonsimple irreducible maximal proper subgroup of $G_n$.
\end{enumerate}

Using similar methods these results were improved in \cite{Oz}.
Nevertheless, the cases when $G$ is a simple Lie group and $\rho:G\to G_n$
is an irreducible representation have remained unanswered.
As a consequence of our main result we will show that if $G$ is
one of the classical Lie groups $SO(k)$, $SU(k)$, $Sp(k)$ such a reduction is 
impossible except for the case $\SU(3)\hookrightarrow \SO(6)$ 
(and the obvious case $G_n=G$).

The paper is organized as follows.
The main results are described in the next section. Homotopy theoretical
results needed in their proofs are contained in Section 3. The proofs
themselves appear in Section 4, based on statements
on dimensions of real representations of classical Lie groups.
The computations proving these statements are carried out in the following
section. In an appendix we give precise self-contained proofs of several
more or less known results which we need and which may be of independent 
interest; we also discuss possible generalizations of our main results.

\section{Main results}

To state our main result we recall several definitions and theorems.
For any prime $p$ let $\nu_p$ stand for the $p$-adic valuation.

The Hurwitz-Radon number $a(r)$ is the power of $2$ given by
$$
\nu_2(a(r)) =
\#\{ i \mid 1\le i\le r-1\,\,\ \text{and $i\equiv 0,1,2,4$ mod 8}\}. 
$$
For $n+1=(2l+1)2^{\beta+4\gamma}$ with $\beta\in\{0,1,2,3\}$  put
$j(n)=2^\beta+8\gamma$. A well known result of Adams 
(\cite{Adams1}) says that the structure group $\SO(n)$ in the fibration 
(\ref{pf1}) can be reduced to the standard subgroup $\SO(k)$ if and only if 
$k\ge n-j(n)+1$.
This condition can be expressed in terms of the Hurwitz-Radon numbers as
\begin{equation*}n+1\equiv 0 \mod a(n+1-k).
\end{equation*}

The complex James number $b(r)$ is the positive integer with
$$
\nu_p(b(r))=\begin{cases}
\max\{i+\nu_p(i) \mid 1\le i\le (r-1)/(p-1)\} &\text{if $r\ge p$,}\\
0&\text{if $r<p$,}
\end{cases}
$$
for all primes $p$.

Similarly, the quaternionic James number $c(r)$ is 
the positive integer determined by
\begin{align*}
\nu_2(c(r))&=
\max\{2r-1,\, 2i+\nu_2(i) \mid 1\le i\le r-1\}\, ,\\
\nu_p(c(r))&=\nu_p(b(2r))\quad\text{for all odd primes $p$}.
\end{align*}

In the papers \cite{AW}, \cite{AT}, \cite{D}, \cite{O1}, \cite{O2} and 
\cite{SS} the problem of
when the structure group $G_n$ in one of the fibrations (1.1) -- (1.3)
can be reduced to a group 
$G=\SU(k)$ or $\Sp(k)$ via a standard inclusion $G\hookrightarrow G_n$ was 
solved and the results were expressed in terms of complex and quaternionic 
James numbers in a way similar to the result quoted above.
The following theorem  can be regarded as a generalization of these results. 

\begin{thm}\label{main}
Let $G_n$ be one of the groups $\SO(n)$, $\SU(n)$ or $\Sp(n)$ and let
$G$ be one of the groups $\SO(k)$ with $k\ge 4$,
$\SU(k)$ with $k\ge 2$ or $\Sp(k)$ with $k\ge 2$. Suppose that the dimension of 
the sphere $G_{n+1}/G_n$ is at least $8$ and that $\dim G<\dim G_n$.
Then the structure group $G_n$ of the principal fibre bundle 
\begin{equation}\label{pf}
G_n\to G_{n+1}\to G_{n+1}/G_n
\end{equation}
can be reduced to $G$ via a homomorphism $\rho:G\to G_n$
if and only if one of the following cases occurs.
\begin{enumerate}
\item[(A)] $G_n=\SO(n)$, $G=\SO(k)$, $n=m-1$, 
$m\equiv 0 \mod a(m-k)$
and, up to conjugation, $\rho$ is the standard inclusion
$\SO(k)\hookrightarrow \SO(n)$.

\item[(B)] $G_n=\SO(n)$, $G=\SU(k)$, $n=2m-1$,
$m\equiv 0 \mod 2^{\nu_2(b(m-k))}$
and, up to conjugation, $\rho$ is the composition of the standard inclusions
$\SU(k)\hookrightarrow \SO(2k)\hookrightarrow \SO(n)$ or the
composition $\SU(4)\to \SO(8)\times \SO(6)\hookrightarrow \SO(15)$
where the first homomorphism is given on the first factor by the standard 
inclusion
and on the second factor by the double covering $\SU(4)\cong \Spin(6)\to 
\SO(6)$.

\item[(C)] $G_n=\SO(n)$, $G=\Sp(k)$, $n=4m-1$,
$m\equiv 0 \mod 2^{\nu_2(c(m-k))}$
and, up to conjugation, $\rho$ is the composition of the standard
inclusions
$\Sp(k)\hookrightarrow SO(4k)\hookrightarrow \SO(n)$
or the exterior square $\Sp (3)\to \SO (15)$.

\item[(D)] $G_n=\SU(n)$, $G=\SU(k)$, $n=m-1$,
$m\equiv 0 \mod b(m-k)$
and, up to conjugation, $\rho$ is the standard inclusion
$\SU(k)\hookrightarrow \SU(n)$.

\item[(E)] $G_n=\SU(n)$, $G=\Sp(k)$, $n=2m-1$,
$m\equiv 0 \mod c(m-k)$
and, up to conjugation, $\rho$ is  the composition of the standard
inclusions
$\Sp(k)\hookrightarrow \SU(2k)\hookrightarrow \SU(n)$.

\item[(F)] $G_n=\Sp(n)$, $G=\Sp(k)$, $n=m-1$,
$m\equiv 0 \mod c(m-k)$
and, up to conjugation, $\rho$ is the standard inclusion
$\Sp(k)\hookrightarrow \Sp(n)$.
\end{enumerate}
\end{thm} 

As a consequence  of Theorem \ref{main} we get a partial answer to Leonard's
question from \cite{L}.

\begin{cor}\label{cor}
Under the assumptions of Theorem {\rm \ref{main}} there is no irreducible 
representation $\rho:G\to G_n$ such that
the structure group $G_n$ of the principal fibre bundle
{\rm (\ref{pf})} can be reduced to $(G,\rho)$.
\end{cor}

\begin{remark}\label{spin} As for $G=\Spin(k)$ or an exceptional simple 
Lie group, using Proposition \ref{lowdim} and similar statements for $\SU(n)$
and $\Sp(n)$ one can easily prove that the 
structure group $G_n$ of (\ref{pf}) cannot be reduced to $G$ via any 
homomorphism $\rho$ with the possible exception of finitely many cases in which 
a certain dimension condition (the same as or similar to that in Proposition 
\ref{lowdim}) is not satisfied.
\end{remark}

\begin{remark}\label{note}
In \cite{BK} and \cite{Th} it was shown that every stably parallelizable 
manifold of dimension $n$ either is parallelizable or has the same span
as $S^n$. This suggests the possibility of extending Theorem \ref{main}
for $G_n=\SO(n)$ from the case of the tangent bundle over $S^n$ to the case 
of a stably trivial, but non-trivial, $n$-dimensional vector bundle over
a stably parallelizable $n$-manifold.
We say a little more about this question in the final section.
\end{remark}

\section{Auxiliary results}

In this section we will summarize the results which will be
needed for the proofs of Theorem \ref{main} and Corollary
\ref{cor} in the next section.
In what follows we will not distinguish between maps and their
homotopy classes.

\begin{prop}\label{lowdim} Let $\tau$ be a principal $\SO(n)$-bundle  
over the suspension $\Sigma X$ of a pointed finite complex $X$.
Suppose that the structure group $\SO(n)$ of $\tau$ can be reduced to 
$(G,\rho )$, where
$\rho:G\to \SO(n)$ is a homomorphism from a Lie group $G$ of dimension less 
than $n-j$, $1\le j<n$.
Then the structure group of $\tau$ can be reduced to the standard
subgroup $\SO(n-j)$ of $\SO(n)$.
\end{prop}

\begin{proof}
We denote the classifying map $X\to \SO(n)$ of a principal fibre bundle $\tau$
by the same letter. 
Suppose that $\SO(n)$ structure can be reduced to $(G,\rho)$. Since
the standard inclusion $\iota:\SO(n-j)\to \SO(n)$ is an
$(n-j-1)$-equivalence and $\dim G\le n-j-1$, the map 
$\rho:G\to \SO(n)$ can be factored as a
composition of a map $\eta:G\to \SO(n-j)$ and $\iota$. Then $\tau$
can be factored through the standard inclusion $\iota$ as
shown by the following diagram.
$$
\xymatrix{
{}&{}&\SO(n-j)\ar@/^/[ddl]^{\iota} \\
X \ar@/^/[urr] \ar[r] \ar[dr]_{\tau}
& G \ar[ur]_{\eta} \ar[d]^{\rho}  & {}\\
{}& \SO(n) &{}
}
$$
\end{proof}

\begin{lemma}\label{lift}
Let $1\le k<n$ and let $n\ge 9$ be odd.
Consider a homomorphism $\rho : G_k \to \O (n)$.
Suppose that there is a map $\rho'' : BG_\infty \to B\O$
such that $\rho''\circ \iota_k\simeq \iota_n\circ B\rho$. A choice of homotopy
induces a diagram of fibrations:
\begin{equation*}\label{diagram}
\xymatrix{
G_\infty /G_k \ar[r]^{\rho'} \ar[d]
& O/O(n) \ar[d]\\
BG_k \ar[r]^{B\rho} \ar[d]_{\iota_k}
& BO(n) \ar[d]^{\iota_n}\\
BG_\infty \ar[r]^{\rho''}
& BO
}
\end{equation*}
Let $\xi\in
\pi_n(BG_k)$ be a homotopy class such that $\tau=B\rho\circ\xi\in
\pi_n(BO(n))$ classifies the fibration {\rm (1.1)}. Then there is
an element $\xi'\in\pi_n(BG_k)$ such that $B\rho\circ\xi'=\tau$ and
$\iota_k\circ\xi'=0$.
\end{lemma}

\begin{proof}
We deal separately with the three cases (a) $G_k=\SO (k)$,
(b) $G_k=\SU (k)$, (c) $G_k=\Sp (k)$.

\par\noindent
{\em Case (a).}\quad
We consider the only non-trivial case: $n\equiv 1$ mod~$8$.
Let $\alpha$ denote the generator of $\pi_n(B\O)= \Zz /2$.  Then
$\pi_{n}(B\O(n))=\Zz /2\oplus\Zz /2$
is generated by $\tau$ and a class $\beta$ such that
$$
\iota_n\beta=\alpha,\quad \iota_n\tau=0
$$
in $\pi_{n}(B\O)$.

We deal first with the case that $k\ge 6$.
Suppose that $\xi\in \pi_{n}(B\SO(k))$ satisfies
$$
B\rho\circ\xi=\tau\quad\text{and}\quad \iota_k\xi=\alpha.
$$

According to \cite{DaMa} (see also the Appendix, Proposition \ref{so6}) 
$\alpha\in \pi_n(B\O)$
can be factored as a composition $\iota_6\eta$ where $\eta\in
\pi_n(B\SO(6))$ and $\iota_6:B\SO(6)\to B\O$ is the standard inclusion. 
This gives the diagram
$$
\xymatrix{
S^{n} \ar@/_/[ddr]_{\eta} \ar@/^/[drr]^{\zeta}
\ar[dr]^{\alpha} & {} & {}\\
{} & B\O & B\SO (k)  \ar[l]_{\iota_k} \\
{} & B\SO (6) \ar[u]_{\iota_6} \ar[ur]_{\iota} & {}
}
$$
in which $\iota:B\SO (6)\to B\SO (k)$ is the standard inclusion and
$\zeta=\iota\eta$.

If we show that $B\rho\circ\zeta=0\in \pi_{n}(B\O (n))$, then
$\xi'=\xi+\zeta$ will satisfy the required conditions, since
\begin{align*}
B\rho\circ\xi'&=B\rho\circ(\xi+\zeta)=
B\rho\circ\xi+B\rho\circ\zeta=\tau+0=\tau ,\\
\iota_k\xi'&=\iota_k(\xi+\zeta)=\alpha+\alpha=0.
\end{align*}
Write $B\rho\circ \zeta=b\beta+c\tau$ with $b,c\in \Zz /2$. Then
$$
b\alpha=\iota_n(b\beta+c\tau+\tau)=\iota_nB\rho(\zeta+\xi)=
\rho''\iota_k(\zeta+\xi)=\rho''(\alpha+\alpha)=0,
$$
which implies that $b=0$. Suppose that $B\rho\circ\zeta=\tau$.
Then $\tau =(B\rho)\circ\iota\circ\eta$ in the diagram:
$$
\xymatrix{
S^{n} \ar@/_/[ddr]_{\eta} \ar@/^/[drr]^{\tau}
\ar[dr]^{\zeta} & {} & {}\\
{} & B\SO (k) \ar[r]_{B\rho} & B\O (n) \\
{} & B\SO (6) \ar[u]_{\iota} &
}
$$
Consequently, the structure group $\O (n)$ of $\tau$ can be
reduced to $\SO (6)$. Since
$$
\dim \SO (6)=15<n-j(n),\quad \text{for }n\equiv 1 \mod 8,\ n\ge 17,
$$
it follows from Proposition \ref{lowdim} that $n=9$. However, the only
possible homomorphism $\SO (6)\to \O(9)$ is (up to conjugation) the
standard inclusion. 
In this case $\rho''=\operatorname{id}$,
which leads to the contradiction:
$$
\alpha=\rho''\iota_6\eta=\rho''\iota_k\iota\eta=\iota_n (B\rho ) \iota\eta
=\iota_n\tau=0.
$$

Finally, if $k<6$, then according to \cite{DaMa}, $(\iota_k)_*:
\pi_{n}(BO(k))\to \pi_{n}(BO)$ always vanishes and we may take
$\xi'=\xi$.

\par\noindent
{\em Case (b).}\quad
In the complex case the assertion is trivial, because $\pi_n(B\SU )=0$.

\par\noindent
{\em Case (c).}\quad
The only non-trivial case occurs for $n\equiv 5$ mod~$8$.
In the Appendix, Proposition \ref{sp1}  we show that the generator of 
$\pi_{8k+4}(\Sp)=\Zz /2$ lifts to an element in
$\pi_{8k+4}(\Sp(1))$. Using this fact the proof proceeds as in (a).
\end{proof}

In the proof of Theorem \ref{main} we will need some properties of stunted
projective and quasiprojective spaces. We will describe them for the real, 
complex 
and quaternionic cases together. Denote by $\Ff$ one of the fields 
$\Rr$, $\Cc$ and $\Hh$ and put $d=1$, $2$ or $4$, respectively, for its 
real dimension. $G(\Ff^n)$ will stand for the group 
$\O(n)$, $\U(n)$ 
or $\Sp(n)$ according to the chosen field $\Ff$. Let us recall that
$G_n$ denotes one of the groups $SO(n)$, $SU(n)$ or $Sp(n)$ and that
\begin{equation*}
G_n/G_k=G(\Ff^n)/G(\Ff^k)\quad \text{for } 1\le k\le n.
\end{equation*}
$P_n$ will stand for the 
projective space $P(\Rr^n)$, $P(\Cc^n)$
or $P(\Hh^n)$ and $Q_n$ for the corresponding quasiprojective space
$Q(\Rr^n)$, $Q(\Cc^n)$ or $Q(\Hh^n)$.
This space is a Thom space $P_n^{\zeta}$ where $\zeta$
is a certain real vector bundle over $P_n$ of dimension $d-1$. 
Denote the Hopf bundle over $P_r$ by $H$ and write $\Ff$ for trivial real 
vector bundle with fibre $\Ff$. The real tangent 
bundle to $P_r$ has the property
$$
\tau(P_r)\oplus\zeta\oplus\Rr = rH^*
$$
where $H^*=\operatorname{Hom}_{\Ff}(H,\Ff)$ is the dual vector bundle to $H$.
Let $t(r)$ denote the order of the bundle $H-\Ff$ in 
$\widetilde J(P_r)$, that is, $t(r)$ is the least integer $\ge 1$ 
such that the sphere bundle of $t(r)H$ is stably fibre homotopy trivial. 
Classical computations in \cite{Adams1}, \cite{AW},
\cite{AT} and \cite{SS} showed that $t(r)=a(r)$, $b(r)$ and $c(r)$ in 
the real, complex and quaternionic cases, respectively.
 
Considering the reflection maps 
\begin{equation*}
\phi: Q_{\infty}/Q_k\to G(\Ff^{\infty})/G(\Ff^k)
\end{equation*} 
and writing $a_i$, $b_i$, $c_i$ for generators of
$\widetilde H_i(Q(\Rr^{\infty});\Zz /2)=(\Zz /2)a_i$, 
$\widetilde H_{2i+1}(Q(\Cc^{\infty}); \Zz )$ $=\Zz b_i$, 
$\widetilde H_{4i+3}(Q(\Hh^{\infty});\Zz)=\Zz c_i$, respectively,
we have inclusions (for unreduced homology)
\begin{multline*}
{\phi}_*: H_*(Q(\Rr^{\infty})/Q(\Rr^{k});\Zz /2)=\Zz/2\oplus
\bigoplus_{i=k}^{\infty}\Zz /2\ a_i\longrightarrow
H_*(O/O(k);\Zz /2)=\\
\Zz /2[a_0,a_1,a_2,\dots]/(a_0=1,a_1=0,\cdots,
a_{k-1}=0,a_k^2=0,
\cdots,a_{i}^2=0,\cdots)
\end{multline*}
where $H_*(O/O(k);\Zz /2)$ is a module over 
$$
H_*(O;\Zz /2)
=\Zz /2[a_0,a_1,a_2,\dots]/(a_0^2=1,a_1^2=0,\cdots),
$$
and for
$\Ff=\Cc$ or $\Hh$
\begin{multline*}
{\phi}_*:H_*(Q(\Ff^{\infty})/Q(\Ff^{k});\Zz)=\Zz \oplus
\bigoplus_{i=k}^{\infty}\Zz\ e_{i}\longrightarrow
H_*(G(\Ff^{\infty})/G(\Ff^k);\Zz )=\\
\Zz [e_{0},e_{1},e_{2},\dots]/(e_{0}=0,\cdots,
e_{k-1}=0,e_k^2=0, \cdots,e_{i}^2=0,\cdots)
\end{multline*}
where $e_{i}$ stands for $b_{i}$ or $c_{i}$, respectively, and 
$H_*(G(\Ff^{\infty})/G(\Ff^k);\Zz )$ is a module over 
$$
H_*(G(\Ff^{\infty});\Zz )=\Zz [e_{0},e_{1},\dots]/(e_{0}^2=0,\cdots).
$$

Since $\SO=\O/\O(1)$ and $\SU=\U/\U(1)$, the formulas above describe also
the homology of $\SO$ and $\SU$.

\begin{prop}\label{hurewicz}
Let $n$ be odd and $k\ge 1$. Consider a map
$f:G_{\infty}/G_k\to O/O(n)$ which fits into a commutative diagram
$$
\xymatrix{
G_{\infty} \ar[d]_{} \ar[r]^{\tilde f}
& O \ar[d]^{} \\
G_{\infty} /G_k \ar[r]^{f} & O/O(n)
}
$$
in which $\tilde f$ is an $H$-map.
Suppose that
\begin{equation}\label{map}
f_*:\pi_n(G_{\infty}/G_{k})\to \pi_n(O/O(n)) = \Zz /2
\end{equation}
is onto. Then $n+1$ is divisible by $d$,
say $n=dm-1$, and $\nu_2(m) \ge \nu_2(t(m-k))$.

The same is true if we replace $G_{\infty}$ and $G_{k}$ by
$G(\Ff^{\infty})$ and $G(\Ff^{k})$.
\end{prop}

\begin{proof}  The homomorphism
$$
f_* : H_n(G_\infty /G_k;\Zz /2) \to H_n(\O /\O (n);\Zz /2)
=(\Zz /2)a_n
$$
maps decomposable elements to $0$, since $f_*$ lifts to a ring homomorphism.
So it follows at once that $n=dm-1$ and that there is an element
$x\in\pi_n(G_\infty /G_k)$ whose Hurewicz image is 
equal to 
$a_{m-1}$, $b_{m-1}$ or $c_{m-1}$ modulo $2$ and products.
Hence the projection
$$
G_\infty /G_k \to G_\infty /G_{m-1}
$$
maps $x$ to a generator of the group $\pi_n(G_\infty /G_{m-1})_{(2)}$,
which is $\Zz /2$ in the orthogonal case, $\Zz_{(2)}$ in the unitary
and symplectic cases.
(The lower index $(2)$ means localization at the prime $2$.)

Now recall that there are stable maps
$$
\theta : G_\infty /G_k \to Q_{\infty}/Q_k
$$
splitting the reflection maps
$\phi : Q_\infty /Q_k \to G_\infty /G_k$.
(See, for example, \cite{Crabb1} for the construction of $\theta$ and
\cite{James} for a description of $\phi$.)
These maps are compatible with the projections $G_\infty /G_k \to
G_\infty /G_l$ and $Q_\infty /Q_k \to Q_\infty /Q_l$ for $k\le l$.

We shall use the letter $\omega$ for stable homotopy:
the symbols $\widetilde\omega_i$ and $\widetilde\omega^i$ will
stand for reduced stable homotopy and cohomotopy groups,
respectively.

Thus $\theta (x)$ gives a class in $\widetilde\omega_n(Q_\infty /Q_k)$
that maps to an odd multiple of the generator of 
$\widetilde\omega_n(Q_\infty /Q_{m-1}) = \Zz /2$, $\Zz$ or $\Zz$, in
the three cases.

The remainder of the proof is an essentially classical computation
using the stable Adams operation $\psi^3$ in $2$-local real $K$-theory,
$KO$.
One may either proceed directly (as we shall show below for the orthogonal
case) or dualize as follows.

By connectivity, the map 
$\widetilde\omega_n(Q_m/Q_k)\to\widetilde\omega_n(Q_\infty /Q_k)$
is surjective. 
So there is an element $y\in \widetilde\omega_n(Q_m/Q_k)_{(2)}$ that
maps to a generator of $\widetilde\omega_n(Q_m/Q_{m-1})_{(2)}=\Zz_{(2)}$.

Now $Q_m/Q_k$ is the Thom space of the bundle $kH^*\oplus\zeta$ over
$P_{m-k} = P(\Ff^{m-k})$.
Its stable dual is, according to \cite{At}, the Thom space of the
virtual bundle $\Rr -mH^*$ over the same space $P_{m-k}$. 
It follows that the restriction map
$$
\widetilde\omega^0(P_{m-k}^{m(\Ff -H^*)})_{(2)} \to 
\widetilde\omega^0(P_{m-(m-1)}^{m(\Ff -H^*)})_{(2)}=
\widetilde\omega^0(S^0)_{(2)}
=\Zz_{(2)}
$$
is surjective.
The Hurewicz image of $y$ in $K$-theory gives under duality a class in
$$
\widetilde KO^0(P_{m-k}^{m(\Ff -H^*)})_{(2)}
$$
that is fixed by $\psi^3$ and restricts to a generator of
$\widetilde KO^0(S^0)_{(2)}=\Zz_{(2)}$.

Alternatively, this shows that the vector bundle $mH^*$ over
$P_{m-k}$ is stably fibre homotopy trivial at the prime $2$
(see \cite{At}, the proof of Proposition 2.8 and \cite{O1}, the proof of 
Theorem 2.2) 
and leads to the equivalent condition that 
$m(\Ff -H^*)\in KO^0(P_{m-k})_{(2)}$ lies in the image of $\psi^3-1$.
(See, for example, \cite{CK3}, Theorem 5.1 modified to $KO$.) 

The proof is completed by calculations of the Adams operation
$\psi^3$. For the unitary and symplectic cases
we refer to \cite{AW} and \cite{SS}.
For the orthogonal case we outline a proof below, which 
is perhaps more direct than the classical calculation.
\end{proof}

We write $kO$ for connective real K-theory.

\begin{prop}\label{ko}
Let $n$ be odd and $1\le k\le n$. 
Suppose that there is an element
$y\in \widetilde kO_n(P(\Rr^\infty )/P(\Rr^k))_{(2)}$ 
which is fixed by the Adams operation $\psi^3$ and maps
to the generator of $\widetilde kO_n(P(\Rr^\infty )/P(\Rr^n))_{(2)}=\Zz /2$.
Then 
$$
n+1 \equiv 0 \mod a(n+1-k)
$$
or, equivalently, $k\ge n-j(n)+1$.
\end{prop}

\begin{proof} 
We deal first with the case that $n+1$ is divisible by $8$.
Using connectivity and duality we may make the identifications
\begin{multline*}
\widetilde KO^0(P(\Rr^{n-k+2})^{(n+1)\Rr - (n+2)H})
= \widetilde KO_n(P(\Rr^{n+2})/P(\Rr^k))\\
=\widetilde kO_n(P(\Rr^{n+2})/P(\Rr^k))
=\widetilde kO_n(P(\Rr^\infty )/P(\Rr^k)).
\end{multline*}
The first group is isomorphic, by Bott periodicity, to
$\widetilde KO^0(P(\Rr^{n-k+2})^{-H})$.
Standard computations of the $K$-groups of real projective spaces 
give that
$\widetilde kO_n(P(\Rr^\infty )/P(\Rr^k))_{(2)}$
is cyclic of order $2a(n+1-k)$ and that $\psi^3$ acts as multiplication
by $3^{(n+1)/2}$.

A generator is fixed by $\psi^3$ if and only if
$$
3^{(n+1)/2}-1 \equiv 0 \mod 2a(n+1-k).
$$
Since $\nu_2(3^{(n+1)/2}-1) = \nu_2(n+1) +1$, the result
follows in this case.

The case $n+1\equiv 4 \mod 8$ is similar. We have
$\widetilde kO_n(P(\Rr^\infty )/P(\Rr^{n-4}))_{(2)}=\Zz /16$ and
$\psi^3$ acts as $3^{(n+1)/2}$. So $k>n-4$.

Finally, for $n+1\equiv 2\mod 4$, we have 
$\widetilde kO_n(P(\Rr^\infty )/P(\Rr^{n-2}))_{(2)}=\Zz /2$ and
the projection map to
$\widetilde kO_n(P(\Rr^\infty )/P(\Rr^n))_{(2)}=\Zz /2$ is zero.
So $k>n-2$.

(It is traditional to use mod $2$ homology and Steenrod operations
for the last two steps, but $kO$-theory provides a uniform proof.)
\end{proof}

Let $\Kk =\Cc$ or $\Hh$.
Suppose that  the field $\Ff =\Cc$ or $\Hh$ is a vector space 
over $\Kk$.
Put again $t(r)=b(r)$ for 
$\Ff=\Cc$ and $t(r)=c(r)$ if $\Ff=\Hh$. Let $n$ be odd
for $\Kk=\Cc$ and $n\equiv 3 \mod 4$ for $\Kk=\Hh$.
The following statement may be established by the same method
as Proposition \ref{hurewicz}.

\begin{prop}\label{othercases} 
Consider a map  $f:G(\Ff^{\infty})/G(\Ff^k)\to 
G(\Kk^{\infty})/G(\Kk^n)$ which lifts to an H-map
$\tilde f : G(\Ff^{\infty})\to G(\Kk^{\infty})$. If
\begin{equation*}
f_*:\pi_n(G(\Ff^{\infty})/G(\Ff^k))\to 
\pi_n(G(\Kk^{\infty})/G(\Kk^n))\cong \Zz ,
\end{equation*}
is onto, then $n=m\operatorname{dim}_{\Kk}{\Ff}-1$ 
and $m\equiv 0\mod t(m-k)$.
\end{prop}

\begin{remark}\label{sufficient} If $f$ in the statement of 
Proposition \ref{hurewicz} 
is the standard inclusion then the condition on $n=dm-1$ is not only necessary 
but also sufficient for 
$f_*:\pi_n(G_{\infty}/G_{k})\to \pi_n(\O/\O(n))$
to be onto. This can be shown by reversing the proof, since the 
conditions on $n$ and $k$ ensure that we are in stable range. The same 
applies to
Proposition \ref{othercases}. 
\end{remark}

\section{Proofs}

Let $\lambda^i$ stand for the representation given by the $i$-th
exterior power and let $\bar{\ }$ denote complex
conjugation. 

\begin{proof}[Proof of Theorem~{\rm\ref{main}}]
According to \cite{L} the structure group $G_n$ in (\ref{pf}) cannot be 
reduced to any proper subgroup for even $n\ge 8$. For odd $n\ge 9$ we will  
examine different $G_n$ and $G$ separately.

\medskip
\noindent {\bf A.} Let $k\le n$. Consider a homomorphism 
$\rho:\SO(k)\to \SO(n)$ which reduces the structure group $\SO(n)$ of the fibre
bundle (\ref{pf1}). First, suppose that the class of the representation 
$\rho$ in $RO(\SO(k))$ is a polynomial in exterior powers. 

\begin{lemma}\label{ext}
Let $\rho: G=\SO (k)\to \O(n)$ be a homomorphism which extends to
a homomorphism $\O (k)\to \O (n)$.
Then there is a map $\rho'':B\SO \to B\O$ such that $\rho''\circ \iota_k
\simeq \iota_n\circ B\rho$, where $\iota_k : B\SO (k)\to B\SO$ and
$\iota_n : B\O (n)\to B\O$ are the standard inclusions.
\end{lemma}

\begin{proof} The real representation ring $RO(\O(k))$ is 
generated by the exterior powers $\lambda^i$ of the basic representation.
Given $\rho$ there is a polynomial $p$ in exterior powers such that
$$p(\lambda^1,\lambda^2,\dots,\lambda^k)(\xi-k)=\rho(\xi)-n$$
for any vector bundle $\xi$ of dimension $k$. This polynomial defines $\rho''$.
To show the existence of such a polynomial it is enough to consider the 
special case
$\rho=\lambda^i$, for which we may take $p$ to be equal to $\sum_{j=0}^{i-1}
\binom{k}{j}\lambda^{i-j}$. 
\end{proof}

The situation is described by the commutative diagram (of fibres)
\begin{equation*}\label{diagram}
\xymatrix{
\SO /\SO (k) \ar[r]^{\rho'} \ar[d]
& \O/\O(n) \ar[d]\\
B\SO (k) \ar[r]^{B\rho} \ar[d]_{\iota_k}
& B\O (n) \ar[d]^{\iota_n}\\
B\SO \ar[r]^{\rho''}
& B\O
}
\end{equation*}

Suppose that the classifying map $\tau:S^n\to B\O (n)$ for the
fibration (\ref{pf1}) can be written as a composition $\tau=B\rho\circ\xi$
with $\xi:S^n\to B\SO (k)$. According to Lemma \ref{lift} we can
suppose that $\iota_k\circ\xi$ is zero in $\pi_n(B\SO )$. Hence both
$\tau$ and $\xi$ can be lifted to a non-trivial element 
$t\in\pi_n(\O /\O (n))\cong \Zz /2$
and $x\in\pi_n(\SO /\SO (k))$, respectively, such that $t=\rho'\circ
x$. 
Hence the map $\rho'$ satisfies
the assumptions of Proposition \ref{hurewicz}. Consequently, 
$n+1\equiv 0 \mod a(n+1-k)$,
which is equivalent to
$$
k\ge n+1-j(n).
$$
Since $n\ge 9$, the inequality above yields $k\ge 8$.

The homomorphism $\rho$ is a sum of irreducible representations.
If it were different from the standard inclusion, then,
by the Weyl Dimension
Formula (see (i) of Proposition \ref{est1} and Remark \ref{est2}
in the next section) and the inequality above, its dimension 
would be at least 
$$
\min\{\dim 2\lambda^1,\, \dim\lambda^2\}=2k\ge 2(n-j(n)+1)>n,
\quad\text{for all } n\ge 9,
$$
which is a contradiction.

Now suppose that the class of $\rho:\SO(k)\to \SO(n)$ in $RO(\SO(k))$ is not 
a polynomial in exterior powers. In this case we can use the following 
lemma. Its proof is based on the Weyl Dimension Formula and is postponed 
to the next section. 

\begin{lemma}\label{dimso}
Let $k\ge 5$. If the class of a representation $\rho:\SO(k)\to \SO(m)$ 
in $RO(\SO(k))$ is not a polynomial in exterior powers, then for all
$n\ge m$
$$
\dim \SO(k)<n-j(n).
$$
\end{lemma}

According to this lemma and Proposition \ref{lowdim}
the reduction of $\SO(n)$ in (\ref{pf1}) to $(\SO(k),\rho)$ is impossible.
\vskip 2mm

\noindent {\bf B1.} Let $n\ge 9$ be odd, $k\ge 2$. 
\begin{lemma}\label{ext-complex}
Let $\rho: G=\SU (k)\to \O(n)$ be a homomorphism 
whose class in $R(\SU (k))$ (after complexification) 
is of the form $q(\lambda^1,\dots,\lambda^{[k/2]},\overline{\lambda^1},
\dots,\overline{\lambda^{[k/2]}})$ where $q$ is a polynomial such that
\begin{equation}\label{polynomial}
q(x_1,\dots,x_{[k/2]},y_1,\dots,y_{[k/2]})=
q(y_1,\dots,y_{[k/2]},x_1,\dots,x_{[k/2]}).
\end{equation}
Then there is a map $\rho'':B\SU \to B\O$ such that $\rho''\circ \iota_k
\simeq \iota_n\circ B\rho$, where $\iota_k : B\SU (k)\to B\U$ and
$\iota_n : B\O (n)\to B\O$ are the standard inclusions.
\end{lemma}

\begin{proof} 
If the complexification of $\rho:SU(k)\to U(n)$
has the form described above, then an extension $BSU\to BU$ can be constructed 
as in the proof of Lemma \ref{ext} using a polynomial $p$ of the same form. 
This extension can be  
factored through the complexification $BSO\to BSU$.
The reason is that any complex vector bundle of the form 
$\eta\oplus\overline\eta$ is the 
complexification of the realification of $\eta$ and any complex vector bundle 
of the form $\eta\otimes\overline\eta$ is the complexification of the real Lie 
algebra bundle of skew-adjoint endomorphisms of $\eta$. 
\end{proof}

Consider a representation $\rho:\SU(k)\to \O(n)$ 
which reduces the structure group $\SO(n)$ of (\ref{pf1}). 

Suppose first that $\rho$ is of the type described in Lemma \ref{ext-complex} 
so that we have a commutative diagram:
\begin{equation*}
\xymatrix{
\SU/\SU(k) \ar[r]^{\rho'} \ar[d]
& \O/\O(n) \ar[d]\\
B\SU (k) \ar[r]^{B\rho} \ar[d]_{\iota_k}
& B\O(n) \ar[d]^{\iota_n}\\
B\SU \ar[r]^{\rho''}
& B\O .
}
\end{equation*}

The classifying map $\tau:S^n\to B\O(n)$ for the
fibration (\ref{pf1}) can be written as a composition $\tau=\rho\circ\xi$
with $\xi:S^n\to B\SU (k)$. By Lemma \ref{lift}
$\tau$ and $\xi$ can be lifted to $t\in\pi_n(\O/\O(n))$
and $x'\in\pi_n(\SU /\SU (k))$, respectively, such that $t=\rho'\circ
x'$ and $t$ is the generator of $\pi_n(\O/\O(n))$. 
So the map $\rho'$ satisfies
the assumptions of Proposition \ref{hurewicz}. Consequently, $n=2m-1$ and 
$m\equiv 0 \mod 2^{\nu_2(b(m-k))}$. For the maximal integer $k$ 
satisfying this condition put $j_2(n)=n+1-2k$. 
Now the divisibility condition above is equivalent to the inequality
$$
2k\ge n+1-j_2(n).
$$
It is clear, by comparing the real and complex lifting problems,
that $j_2(n)\le j(n)$.
Since $n\ge 9$, the inequality above yields $k\ge 4$.

Suppose that $k\ge 5$.
If $\rho$, which is a sum of irreducible representations,
were different from the standard inclusion, then, by the Weyl Dimension
Formula (see (ii) of Proposition \ref{est1} and Remark \ref{est2}) 
and the inequality above, its dimension
would be at least 20 and greater than or equal to
$$
\min\{2\dim_{\Cc}(\lambda^1+\overline{\lambda^1}),
\dim_{\Cc}(\lambda^2+\overline{\lambda^2})\}=
\min\{4k,k^2-k\}=4k\ge 2(n-j(n)+1)>n.
$$
This means that $\rho$ has to be the standard inclusion in this case.

Now consider the case: $k=4$.
In the same way we get
$$
8=2\cdot4\ge n-j_2(n)+1.
$$
Consequently, $n=9$, $11$ or $15$ in this case. 
This allows only two possibilities for $\rho$:
$\lambda^1+\overline{\lambda^1}$ and 
$\lambda^1\cdot\overline{\lambda^1}-1$  of dimension $8$ and $15$, 
respectively. 

The latter homomorphism can be factored via a double covering as
$$
\SU(4)\cong \Spin(6)\to \SO(6)\xrightarrow{\ \lambda^2\ } \SO(15).
$$ 
However, the reduction of $\SO(15)$ to $\SO(6)$ was excluded in {\bf A}. 

Hence $\rho$ is a standard inclusion corresponding to 
$\lambda^1+\overline{\lambda^1}$.

\vskip 2mm

\noindent {\bf B2.}
Consider a homomorphism $\rho:\SU(k)\to \O(n)$ whose class in $R(\SU(k))$ is
not of the form described by (\ref{polynomial}). In the next section we prove

\begin{lemma}\label{dimsu}
Let $k\ge 2$. If the class of an irreducible representation $\SU(k)\to \O (m)$
in $R(\SU(k))$ is not  a polynomial in exterior powers of the form
{\rm(\ref{polynomial})}, then for all
$n\ge m$
$$
\dim \SU(k)<n-j(n)
$$
with the just two exceptions: the double covering 
$\SU(4)\cong \Spin(6)\to \SO(6)$
which after complexification gives $\lambda^2:\SU(4)\to \SU(6)$, and
the representation $\SU(8)\to \SO(70)$ which after complexification gives
$\lambda^4:\SU(8)\to \SU(70)$.
\end{lemma}

According to Proposition \ref{lowdim} this lemma excludes reductions 
to $(\SU(k),\rho)$ apart from two possible exceptional cases:
that after complexification one of the irreducible summands of $\rho$ is  
$\lambda^2:\SU(4)\to \U (6)$ and $n\le 23$ (for $n\ge 25$ we can use
Proposition \ref{lowdim}) or  
$\lambda^4:\SU(8)\to \U (70)$ with $n=71$ (for $n\ge 73$ we can again use
Proposition \ref{lowdim}).   

The representation $\rho$  of $\SU(4)$ containing $\lambda^2$ as a summand 
has to contain as a summand also 
$\lambda^1+\overline{\lambda^1}$, otherwise it could be factored through 
$\SO(6)$, which is impossible. So the representation $\rho$ of $\SU(4)$
has to be of the form
\begin{equation}\label{su4}
\SU(4)\to \SO(8)\times  \SO(6) \to 
\SO(14)\times \SO(n-14)\hookrightarrow \O(n)
\end{equation}
where the first homomorphism is given by the standard inclusion
and the double covering and the second one contains the standard
inclusion $\SO(8)\times \SO(6)\hookrightarrow \SO(14)$ as the first component. 
However, according to Theorem 2.A  in \cite{L} the structure group $\SO(n)$ 
cannot be reduced to $\SO(14)\times 
\SO(n-14)$ for  $17\le n\le 23$.  (This can be seen by observing that a 
reduction to the product would imply reduction to one of the factors. 
Compare the proof of Lemma \ref{real}.) 

\begin{lemma}\label{real} The
structure group $\SO(15)$ of the principal fibre bundle {\rm(\ref{pf1})} can
be reduced to $\SU(4)$ both through the standard inclusion 
and through the homomorphism of the form {\rm(\ref{su4})}.
\end{lemma}

\begin{proof} By \cite{D} or by {\bf B1} and Remark \ref{sufficient} the 
structure group $\SO(15)$ of (\ref{pf})
can be reduced to $\SU(4)$ using the standard inclusion.
Since the tangent bundle to $S^{15}$ is determined by the non-trivial
element in $\pi_{14}\big(\SO(15)\big)\cong \Zz /2$, there is
a map $\gamma : S^{14}\to \SU(4)$ such that the composition with
$\SU(4)\hookrightarrow \SO(15)$ is non-trivial. The composition of
$\gamma$ with any homomorphism
$$
\SU(4)\to \SO(6)\hookrightarrow \SO(15)
$$
is trivial since the inclusion  $\SO(6)\hookrightarrow \SO(15)$
induces the trivial homomorphism
$$
\pi_{14} \big(\SO(6)\big) \to \pi_{14}\big(\SO(15)\big)
$$
by {\bf A}.

Hence the composition
$$
S^{14} \xrightarrow{\ \gamma\ } \SU(4)\hookrightarrow \SO(8)\times
\SO(6) \to \SO(14)\hookrightarrow \SO(15)$$
is homotopic to the product (in $\SO(15)$) of two maps
$S^{14}\to \SO(15)$, one of which is non-trivial and the other trivial.
That is why this composition is non-trivial. So $\gamma :
S^{14}\to  \SU(4)$ determines the reduction  through the homomorphism
(\ref{su4}).
\end{proof}

\begin{lemma}\label{su8}
Let $n=71$. Then the structure group $\SO(71)$ in {\rm (\ref{pf1})}
cannot be reduced
to $\SU(8)$ via the homomorphism induced by $\lambda^4:\SU(8)\to \SU(70)$.
\end{lemma}

\begin{proof} 
Suppose that a reduction does exist. 
Then the classifying map for (\ref{pf1}) admits a factorization 
$$
S^{70}\xrightarrow{f} \SU(8)\xrightarrow{\rho} \SO(70)\hookrightarrow \SO(71).
$$
Denote by $\eta$ the real vector bundle 
over $\Sigma \SU(8)$ which is classified by
$\rho$. Then the tangent bundle to $S^{71}$ is isomorphic to 
$f^*(\eta)\oplus\Rr$. To prove that the factorization above is 
impossible it 
is sufficient to show that the Stiefel-Whitney class $w_{64}(\eta)=0$.
Since $\dim \Sigma \SU(8)=64$, in this case the vector bundle $\eta$ would have
$7$ linearly independent vector fields, and, consequently, there would be $8$ 
linearly independent vector fields on $S^{71}$, which is a contradiction.

Now $\eta$ is the pullback, via the classifying map
$\Sigma \SU(8)\to B\SU(8)$, of a bundle $\widetilde\eta$
over $B\SU(8)$. The cohomology ring $H^*(B\SU(8);\, {\Zz}/2 )
= \Zz/2[c_2,c_3,\ldots ,c_8]$ is polynomial
on the mod $2$ reductions of the Chern classes of the universal bundle.
As the dimension of the generators is bounded by $16$,
$w_{64}(\widetilde\eta)$ is a product of lower dimensional classes.
Since the products in $H^*(\Sigma \SU(8);\Zz)$ are trivial,
it follows that $w_{64}(\eta )=0$.
\end{proof}

\vskip 2mm

\noindent {\bf C.} Consider a representation $\rho:\Sp(k)\to \O(n)$
which reduces the structure group $\SO(n)$ of (\ref{pf1}).

\begin{lemma}\label{ext-quaternionic}
Let $\rho: G=\Sp (k)\to \O(n)$ be a homomorphism.
Then there is a map $\rho'':B\Sp \to B\O$ such that $\rho''\circ \iota_k
\simeq \iota\circ B\rho$.
\end{lemma}

\begin{proof} 
In $R(\Sp(k))$ any virtual representation is self-conjugate and 
corresponds to a polynomial in exterior powers:  $\lambda^i$ is real if 
$i$ is even
and quaternionic if $i$  is odd. Real virtual 
representations in $RO(\Sp(k))$ are given by those polynomials in $R(\Sp(k))$
which have even coefficients at monomials 
$(\lambda^1)^{j_1},(\lambda^2)^{j_2},\dots,(\lambda^k)^{j_k}$  with
$j_1+j_2+\dots+j_k\equiv 1$ mod 2. The proof can be completed as in the proof 
of Lemma \ref{ext-complex}.
\end{proof}

The classifying map
$\tau:S^n\to BO(n)$ for the
fibration (\ref{pf1}) can be written as a composition $\tau=\rho\circ\xi$
with $\xi:S^n\to B\Sp(k)$. According to Lemma \ref{lift} this map
can be chosen in such a way that $\iota_k\circ\xi=0\in\pi_n(B\Sp)$.
Hence both
$\tau$ and $\xi$ can be lifted to $t\in\pi_n(\O/\O(n))$
and $x'\in\pi_n(\Sp/\Sp(k))$, respectively, such that $t=\rho'\circ
x'$ and $t$ is the generator of $\pi_n(\O/\O(n))$. 
So the map $\rho'$ satisfies
the assumptions of Proposition \ref{hurewicz}. Consequently, $n=4m-1$ and 
$m\equiv 0 \mod 2^{\nu_2(c(m-k))}$. For maximal $k$ satisfying this condition
put $j_4(n)=n+1-4k$. Now the divisibility condition above is equivalent to
the inequality
$$
4k\ge n+1-j_4(n).
$$
Again it is clear that $j_4(n)\le j_2(n)$.
Since $n\ge 9$, the inequality above yields $k\ge 3$.

Consider $k\ge 4$. 
If $\rho$ were not the standard inclusion, then by the Weyl Dimension
Formula (see (iii) of Proposition \ref{est1} and Remark \ref{est2}) 
its dimension would be at least 
$$
\min\{2\dim_{\Cc} 2(\lambda^1),\dim_{\Cc}(\lambda^2 -1)\}=
\min \{8k,(2k+1)(k-1)\}\ge 27.
$$
For $n\ge 27$ the inequality above yields 
$$
(2k+1)(k-1)>2(k-1)^2\ge\frac{(n-j_4(n)-3)^2}{8}\ge n
$$
and
$$
8k\ge 2(n-j_4(n)+1)>n.
$$
This means that $\rho$ has to be the standard inclusion if $k\ge 4$.

For $k=3$ the inequality implies that $n\le 15$.
We have two irreducible representations of $\Sp(3)$ of dimension
$\le 15$: the standard inclusion and the one given by 
 the polynomial $\lambda^2-1$ in $R(\Sp(3))$ of
dimension $14$. 

\begin{lemma}\label{sp3}
Let $n=15$. The structure group $\SO(15)$ in {\rm(\ref{pf1})} can be 
reduced to $\Sp(3)$ via the homomorphism $\rho$ induced by 
$\Sp(3) \xrightarrow{\ \lambda^2\ } \SO(15)$.
\end{lemma}

\begin{proof} Since $\Sp/\Sp(3)$ is 14-connected, $\SO(15)$ can be reduced to
$\Sp(3)$ via $\rho$ if and only if 
$$
\rho'_* : \widetilde H_{15}(\Sp/\Sp(3);\Zz /2)\to 
\widetilde H_{15}(\O/\O(15);\Zz /2)
$$
is onto.
We shall show that the composition
\begin{multline}\label{comp}
\widetilde H_{15}(Q(\Hh^4);\Zz /2) \to \widetilde H_{15}(\Sp (4);\Zz /2)
\to \widetilde H_{15}(\O ;\Zz /2) \\ 
\to \widetilde H_{15} (\O /\O (15);\Zz /2) =\Zz /2
\end{multline}
induced by the reflection map $\phi : Q(\Hh^4) \to \Sp (4)$ and 
$\rho'' : B\Sp (4) \hookrightarrow B\Sp\to B\O$ is non-zero.
On $B\Sp (4)$, $\rho''$ is given on a $4$-dimensional $\Hh$-vector bundle
$\xi$ as the virtual bundle $\lambda^2 \xi - \xi\otimes_\Cc \Hh + \Cc^3$
(with its real structure).

The composition (\ref{comp}) is given by a cohomology class 
$w\in \widetilde H^{15}(Q(\Hh^4);\Zz /2)$. To compute the class it is 
convenient
to lift from $Q(\Hh^4)$, which is the Thom space of the $3$-dimensional
Lie algebra bundle $\zeta$ over the quaternionic projective space
$P(\Hh^4)$, to the sphere bundle $S(\Rr\oplus\zeta )$.
The map
$$
S(\Rr\oplus\zeta) \to P(\Hh^4)^\zeta = Q(\Hh^4) \to \Sp (4)
$$
determines a $4$-dimensional $\Hh$-vector bundle $\xi$ over
$S^1\times S(\Rr\oplus\zeta )$.  
The class $w$ lifts to 
$$
w_{16}(\rho''(\xi )) \in H^{16}(S^1\times S(\Rr\oplus\zeta);\Zz /2).
$$

Now we can write $H^*(S^1;\Zz /2)=\Zz /2[t]/(t^2)$,
$H^*(P(\Hh^4);\Zz /2)=\Zz /2[x]/(x^4)$, where $x=w_4(H)$
is the mod $2$ Euler class of the quaternionic Hopf bundle $H$,
and $H^*(S(\Rr\oplus\zeta);\Zz /2)=\Zz /2[x,y]/(x^4,y^2)$,
where $y$ is the $3$-dimensional class corresponding to the Thom class
of $\zeta$.
The vector bundle $\xi$, constructed using the reflection map,
is a direct sum $\eta\oplus H^\perp$, where $H^\perp$ is the orthogonal
complement of $H$ in the trivial bundle $\Hh^4$ over $P(\Hh^4)$ and
$\eta$ is the quaternionic line bundle obtained by twisting $H$.
It follows that
$$
\rho''(\xi ) =\lambda^2\eta + \eta\otimes_{\Cc} H^\perp +\lambda^2 H^\perp
-\eta\otimes_{\Cc}\Hh -H^\perp\otimes_{\Cc} \Hh +\Cc^3.
$$
Using the triviality of the bundles $\lambda^2\eta$ and $H\oplus H^\perp$,
one obtains
$$
\rho''(\xi ) = (\eta - H)\otimes_{\Cc} \Hh^3 + (H- \eta )\otimes_{\Cc} H 
+\Cc^{15}.
$$
It is understood here that each complex bundle which we have
written down has a real structure. We have to compute the Stiefel-Whitney
class $w_{16}$ of the virtual real vector bundle so defined.
This will be done by calculating the total Stiefel-Whitney classes
of the various constituents: $\eta\otimes\Hh$, $H\otimes\Hh$,
$H\otimes H$ and $\eta\otimes H$.

To compute the Stiefel-Whitney classes involving $\eta$ it is enough to
consider the restrictions to the subspaces $S(\Rr\oplus\zeta)$,
where $\eta$ coincides with $H$, and $S^1\times S^3$,
where $S^3$ is the fibre of $S(\Rr\oplus\zeta)$ 
at a point in $P(\Hh^4)$.
The second restriction gives us the generator of
$\pi_3(\Sp(1))=\pi_4(B\Sp (1))$, which determines a
$4$-dimensional real vector bundle over $S^4$ with non-zero
Stiefel-Whitney class.
(One can, for example, think of this vector bundle as the Hopf bundle
over $S^4=P(\Hh^2)$.) 

One finds that:
$$
w(\eta\otimes\Hh ) = 1+x+ty,\quad
w(H\otimes\Hh ) = 1+x,\quad
w(H\otimes H)=1,\quad
w(\eta\otimes H)=1+ty.
$$
Hence $w(\rho''(\xi )) = (1+x+ty)^3(1+x)^{-3}(1+ty)^{-1}
=1+(x+x^2+x^3)ty$. This verifies that $w_{16}=x^3ty$ is non-zero,
as claimed.
\end{proof}

\vskip 2mm

\noindent {\bf D}, {\bf E}, {\bf F}.
Now suppose $G_n=\SU(n)$ or $\Sp(n)$. Put $m=n+1$ and 
$d=2$ or $4$ in the complex or quaternionic case, respectively. Consider
the diagram
\begin{equation*}
\xymatrix{
G_n \ar[r] \ar[d] & \SO(dm-d) \ar[r]
& \SO(dm-1) \ar[d]\\
G_{n+1} \ar[rr] \ar[d]  &
& \SO(dm) \ar[d]\\
S^{dm-1} \ar[rr]^{=} &
& S^{dm-1}
}
\end{equation*}
If $G_n$ in (\ref{pf}) can be reduced to $G$ via $\rho:G\to G_n$, then
$\SO(dm-1)$ in (\ref{pf1}) can be reduced to $G$ via 
\begin{equation}\label{inclusion}
G\xrightarrow{\rho}G_n\hookrightarrow \SO(dm-1).
\end{equation}
According to the previous steps this composition has to be a standard inclusion
or the composition $\SU(4)\to \SO(8)\times \SO(6)\hookrightarrow \SO(15)$ 
described in {\bf B1} or $\lambda^2:\Sp(3)\to \SO(15)$ from Lemma \ref{sp3}. 
However, the last two homomorphisms cannot be factored 
through $\SU(7)$ or $\Sp(4)$. The inclusion (\ref{inclusion}) has to satisfy 
the 
inequality
\begin{align*}
k&\ge dm-j(dm-1)    & \text{for }G&=\SO(k),\\
2k&\ge dm-j_2(dm-1) & \text{for }G&=\SU(k),\\
4k&\ge dm-j_4(dm-1) & \text{for }G&=\Sp(k),
\end{align*} 
which implies that the cases $G=\SO(k)$ with $G_n=\SU(n)$ or $\Sp(n)$
and $G=\SU(k)$ with $G_n=\Sp(n)$ cannot occur, since $k>n$. In the remaining 
cases
$\rho$ is forced to be a standard inclusion, again for dimensional reasons. 
Then Proposition \ref{othercases} 
gives the divisibility conditions in (D), (E) and (F) of Theorem 2.1.
and Remark \ref{sufficient} says that these conditions are sufficient for
the existence of reductions.
\end{proof}

\section{The Weyl Dimension Formula}

The estimates of dimension of real irreducible representations used in the
previous section are based on the Weyl Dimension Formula. Here we show how
it is used. In particular, we deduce Lemmas \ref{dimso} and \ref{dimsu} from
the following propositions:

\begin{prop}\label{est1}
Let $\rho : G_k\to SO(m)$ be an irreducible, non-trivial, non-standard
real representation.
\smallskip
\item[(i)]
$G_k=SO(k)$, $k\ge 7$.
Then $m\ge k(k-1)/2$.
\item[(ii)]
$G_k=SU(k)$, $k\ge 5$.
Then $m\ge k(k-1)$.
\item[(iii)]
$G_k=Sp(k)$, $k\ge 3$.
Then $m\ge k(2k-1)-1$.
\end{prop}

\begin{remark}\label{est2} The bounds on $m$ are achieved by
(i) $\lambda^2$,
(ii) the underlying real representation of $\lambda^2$,
(iii) the real representation $\lambda^2$ modulo the one-dimensional
trivial summand given by the defining symplectic form.
\end{remark}

\begin{prop}\label{est3} Let $\rho :G_k \to SO(m)$ be a real representation.
\smallskip
\item[(i)] $G_k=SO(k)$.
If $[\rho ]\in RO(G_k)$ is not a polynomial in the exterior powers
$\lambda^j$, $j\geq 1$, then $k\equiv 0$ mod $4$ and
$m\ge \frac{1}{2}\binom{k/2}{k/4}$
$(= \dim \lambda^{k/2}_+)$.
\item[(ii)] $G_k=SU(k)$.
If $[\rho ]\in RO(G_k)$ is not of the type described in 
Lemma {\rm\ref{ext-complex}}
then $k$ is even and $m\ge \binom{k}{k/2}$
$($the dimension of $\lambda^{k/2}$ with its real structure$)$.
\end{prop}

We start by recalling necessary prerequisities from the
representation theory of Lie groups. Let $G$ be a simple Lie
group. Denote by $\omega_1,\omega_2,\dots,\omega_l$ its fundamental
weights. Then any complex irreducible representation is
determined by its dominant weight
$$\omega=m_1\omega_1+m_2\omega_2+\dots+m_l\omega_l,\quad
m_i\in\mathbb N.$$
The complex vector space on which it acts will be
denoted by $V(\omega)=V(m_1,m_2,$ $\dots,m_l)$. Put
$\delta=\sum_{i=1}^{l}\omega_i$.
The Weyl Dimension Formula gives the complex dimension of $V(\omega)$
$$
\dim_{\Cc}V(\omega)=\prod_{\beta>0}\frac{\langle\beta,\omega+\delta\rangle}
{\langle\beta,\delta\rangle},
$$
where $\beta$ goes through all the positive roots. Treating
separately all the types of classical simple Lie groups we always get that
if $ m'_i\ge m_i$ for all $i$, then
\begin{equation}\label{eneq}
\dim_{\Cc}V(m'_1,m'_2,\dots,m'_l)\ge
\dim_{\Cc}V(m_1, m_2,\dots, m_l).
\end{equation}

We will need to compute dimensions of real irreducible
representations. According to Cartan's Theorem (\cite{GG}, page
366) there are two kinds of real irreducible representations:
those the complexification of which is irreducible, and those the
complexification of which is a sum of an irreducible complex
representation with its complex conjugate representation. So the
dimensions of real irreducible representations can be computed
from the knowledge of dimensions of complex irreducible
representations.

\vskip 2mm

\begin{proof}[Proof of\/ {\rm(i)} of Propositions\/ {\rm\ref{est1}} 
and\/ {\rm\ref{est3}}]
First, consider representations of $\Spin(k)$ or
$\SO(k)$ with $k=2l+1\ge 7$. Denote by $\varepsilon_i\pm \varepsilon_j$,
$\pm\varepsilon_i$, $1\le i,j\le l$, $i\ne j$, the roots of the
corresponding Lie algebra $\frak{so}(2l+1)$. Its positive roots
are $\varepsilon_i\pm\varepsilon_j$ for $1\le i<j\le l$, and
$\varepsilon_i$ for $1\le i\le l$. The fundamental weights are
$$\omega_1=\varepsilon_1,\ \omega_2=\varepsilon_1+\varepsilon_2,
\dots,\omega_{l-1}=\varepsilon_1+\dots+\varepsilon_{l-1},
\ \omega_{l}=\frac{1}{2}(\varepsilon_1+\dots+\varepsilon_{l}).$$
A complex representation with dominant weight 
$\omega=\sum_{i=1}^{l}m_i\omega_i$ acting on $V(\omega)$ 
is the complexification of a real representation 
$\rho:\SO(2l+1)\to \SO(m)$ if and only if $m_l$ is even. Such a representation
in $RO(\SO(2l+1))$ is described by a polynomial in exterior powers.
($\omega_j$ corresponds to $\lambda^j$ for $1\le j\le l-1$ and $2\omega_l$
corresponds to $\lambda^l$.)

Put $\omega+\delta=\sum_{i=1}^{l}g_i\varepsilon_i$.
Then the Weyl Dimension Formula reads as
$$\dim_{\Cc} V(\omega)=\prod_{1\le i<j\le l}\frac{(g_i-g_j)(g_i+g_j)}
{(j-i)(2l+1-i-j)}\prod_{1\le i\le l}\frac{2g_i}{2l-2i+1}.$$
From this formula we can see immediately that the inequality
(\ref{eneq}) holds. We compute the dimensions of several irreducible
representations:
\begin{align*}
\dim_{\Cc} V(\omega_{j})&=\binom{2l+1}{j},\quad 1\le j\le l-1,\\
\dim_{\Cc} V(2\omega_{l})&=\binom{2l+1}{l},\quad 
\dim_{\Cc} V(2\omega_{1})=(2l+3)l.
\end{align*}
Let $\rho$ be a real irreducible representation $\SO(2l+1)\to \SO(m)$
with a dominant
weight $\omega$ different from $\omega_1$.
For $l\ge 3$, using the inequality (\ref{eneq}), we get that its dimension
$m=\dim_{\Cc} V(\omega)$ is at least
\begin{equation*}
\min\{\dim_{\Cc} V(2\omega_1),\dim_{\Cc} V(\omega_2),
\dots,\dim_{\Cc} V(\omega_{l-1}),\dim_{\Cc} V(2\omega_l)\}
= (2l+1)l.
\end{equation*}

Now consider representations of $\Spin(k)$ or
$\SO(k)$ with $k=2l\ge 8$. For $l\ge 4$ the roots of the
corresponding Lie algebra $\frak{so}(2l)$ are $\varepsilon_i\pm \varepsilon_j$,
$\pm\varepsilon_i$, $1\le i,j\le l$, $i\ne j$.
The positive roots
are $\varepsilon_i\pm\varepsilon_j$ for $1\le i<j\le l$, and
$\varepsilon_i$ for $1\le i\le l$. The fundamental weights are
\begin{gather*}
\omega_1=\varepsilon_1,\ \omega_2=\varepsilon_1+\varepsilon_2,
\dots,\omega_{l-2}=\varepsilon_1+\dots+\varepsilon_{l-2},\\
\omega_{l-1}=\frac{1}{2}(\varepsilon_1+\dots+\varepsilon_{l-1}-
\varepsilon_{l}),\\
\omega_{l}=\frac{1}{2}(\varepsilon_1+\dots+\varepsilon_{l-1}+
\varepsilon_{l}).
\end{gather*}
A complex representation with dominant weight
$\omega=\sum_{i=1}^{l}m_i\omega_i$ acting on $V(\omega)$
is a  representation
of $\SO(2l)$ if and only if $m_{l-1}+m_l$ is even. 

If $l$ is even, $V(\omega)$ is always 
the complexification of an irreducible representation $\rho:\SO(2l)\to
\SO(m)$. In this case
$RO(\SO(2l))$ is generated by $\lambda^j$, $1\le j\le l-1$, $\lambda^l_+$,
$\lambda^l_-$, where $\lambda^l_-+\lambda^l_+=\lambda^l$.
Here $\omega_j$ is the dominant weight of the complexification of $\lambda^j$ 
for $1\le j\le l-2$,
$\omega_{l-1}+\omega_l$ is the dominant weight of the complexification of 
$\lambda^{l-1}$,
and $2\omega_{l-1}$, $2\omega_{l}$ are the dominant weights 
of complexifications of $\lambda^l_-$, $\lambda^l_+$, respectively.

If $l$ is odd, $RO(\SO(2l))$ is a polynomial ring in exterior powers. 
$\omega_j$ corresponds to the complexification of $\lambda^j$ for $1\le j\le
l-2$ and $\omega_{l-1}+\omega_l$ corresponds to the complexification of
$\lambda^{l-1}$, while the complexification of $\lambda^l$ is the sum of
complex irreducible representations with dominant weights 
$2\omega_{l-1}$ and $2\omega_{l}$.

Let $\omega=\sum_{i=1}^{l}m_i\omega_i$ be the dominant weight of
a complex representation acting on $V(\omega)$. Put $\omega+\delta=
\sum_{i=1}^{l}g_i\varepsilon_i$.
The Weyl Dimension Formula now gives
$$\dim_{\Cc} V(\omega)=\prod_{1\le i<j\le l}\frac{(g_i-g_j)(g_i+g_j)}
{(j-i)(2l-i-j)}.$$
From the formula one readily verifies the inequality (\ref{eneq}) and
computes the dimensions:

\begin{align*}
\dim_{\Cc} V(\omega_{j})&=\binom{2l}{j},\quad 1\le j\le l-2,\\
\dim_{\Cc} V(\omega_{l-1}+\omega_{l})&=\binom{2l}{l-1},\\
\dim_{\Cc} V(2\omega_{l-1})&=\dim_{\Cc} V(2\omega_{l})
=\frac{1}{2}\binom{2l}{l},\\
\dim_{\Cc} V(2\omega_{1})&=(2l-1)(l+1).
\end{align*}
Let $\rho$ be a real irreducible representation $\SO(2l)\to \SO(m)$
the complexification of which has a dominant
weight $\omega$ different from $\omega_1$.
Then according to (\ref{eneq}) for $l\ge 4$ its dimension is either
$\dim_{\Cc} V(\omega)$ or $2\dim_{\Cc} V(\omega)$ and greater than or
equal to
\begin{equation*}
\min\{\dim_{\Cc} V(2\omega_1),\dim_{\Cc} V(\omega_2),
\dots,\dim_{\Cc} V(\omega_{l-1}+\omega_l),
\dim_{\Cc} V(2\omega_{l-1})\}
= l(2l-1).
\end{equation*}

If the class of an irreducible representation $\rho:\SO(k)\to \SO(m)$ is not a
polynomial in exterior powers in $RO(\SO(k))$, then $k=2l\equiv 0$ mod 4
and the complexification of $\rho$ has a dominant weight 
$\omega=\sum_{i=1}^{l}m_i\omega_i$ with $m_{l-1}\ge 2$
or $m_l\ge 2$. By (\ref{eneq}) its dimension is at least 
$$\dim_{\Cc} V(2\omega_{l-1})=\dim_{\Cc} V(2\omega_l)=\frac{1}{2}
\binom{2l}{l}.$$ 
\end{proof}

For $l\ge 4$ and $n\ge \frac{1}{2}\binom{2l}{l}$ we get
$$\dim \SO(2l)=l(2l-1) < n-j(n),$$
which proves Lemma \ref{dimso}.

\begin{proof}[Proof of\/ {\rm(ii)} of Propositions\/ {\rm\ref{est1}} 
and\/ {\rm\ref{est3}}]
Consider $\SU(k)$ with $k\ge 2$.
In standard notation $\varepsilon_i- \varepsilon_j$,
$1\le i,j\le k$, $i\ne j$, are the roots of the
Lie algebra $\frak{su}(k)$. Its positive roots
are $\varepsilon_i-\varepsilon_j$ for $1\le i<j\le k$,
and the fundamental weights are
$$
\omega_1=\varepsilon_1,\ \omega_2=\varepsilon_1+\varepsilon_2,
\dots,\omega_{k-1}=\varepsilon_1+\varepsilon_2+\dots+
\varepsilon_{k-1}.
$$
Let $\omega=\sum_{i=1}^{k-1}m_i\omega_i$ be the dominant weight of
a complex irreducible representation which acts on a complex
vector space $V(\omega)$. The only complex irreducible
representations which are complexifications of
real irreducible representations are those satisfying
$$m_i=m_{k-i},\quad\text{for }1\le i\le k-1,$$
with the exception of the case that $k\equiv 2$ mod 4 and $m_{k/2}$ is odd.
See \cite{S}, Theorem E, page 140.
The remaining complex irreducible representations $V(\omega)$
determine real irreducible representations $\rho$ with
$\dim_{\Rr} \rho=2\dim_{\Cc} V(\omega)$.

The complex exterior power $\lambda^i$ has dominant weight $\omega_i$
and its conjugate representation is $\lambda^{k-i}$. Hence the representations
described in $R(\SU(k))$ by polynomials (\ref{polynomial}) 
are complexifications of real ones. 

Put $\omega+\delta=\sum_{i=1}^{k-1}g_i\varepsilon_i$.
In this case the Weyl Dimension Formula reads as
$$\dim_{\Cc} V(\omega)=\prod_{1\le i<j\le k-1}\frac{(g_i-g_j)}
{(j-i)}\prod_{1\le i\le k-1}\frac{g_i}{k-i}.$$
Again we can check that  the inequality (\ref{eneq})  holds. 
Specific dimensions are:
\begin{align*}
\dim_{\Cc} V(\omega_{j})&=\binom{k}{j},\quad 1\le j\le k-1,\\
\dim_{\Cc} V(2\omega_{1})&=\frac{1}{2}k(k+1)\\
\dim_{\Cc} V(\omega_{1}+\omega_{k-1})&=k^2-1,\\
\dim_{\Cc} V(\omega_{2}+\omega_{k-2})&=\frac{1}{4}k^2(k+1),\quad k\ge 5.
\end{align*}
For $k\ge 5$ consider a real irreducible representation
of $\SU(k)$ which is determined by  a complex irreducible
representation with  a dominant
weight $\omega$ different from $\omega_1$. Its real dimension is at least
\begin{equation*}
\min\{2\dim_{\Cc} V(2\omega_1),2\dim_{\Cc} V(\omega_2),
\dim_{\Cc} V(\omega_1+\omega_{k-1}),
\dim_{\Cc} V(\omega_2+\omega_{k-2})\}=k^2-k.
\end{equation*}

If the class of an irreducible representation $\rho:\SU(k)\to \SO(m)$ (after
complexification) is not a
polynomial in exterior powers in $R(\SU(k))$ of the form (\ref{polynomial}), 
then $k=2l\equiv 0$ mod 4 and the complexification of $\rho$ has 
a dominant weight $\omega=\sum_{i=1}^{k-1}m_i\omega_i$ with $m_{l}\ge 1$. 
According to (\ref{eneq}) its dimension is at least
$$ \dim_{\Cc} V(\omega_{l})=\binom{2l}{l}.$$
\end{proof}

For $l\ge 6$ and $n\ge \binom{2l}{l}$ we get
$$\dim \SU(2l)=4l^2-1 < n-j(n).$$
For $k=4$, the dimension of $\SU(4)$ is 15 and the only representation
of the above form which does not satisfy the required inequality has
the dominant weight $\omega_2$.
For $k=8$, the dimension of $\SU(8)$ is 63 and the only representation
of the above form which does not satisfy the required inequality has
the dominant weight $\omega_4$.
This completes the proof of Lemma \ref{dimsu}.

\begin{proof} [Proof of\/ {\rm(iii)} of Proposition\/ {\rm\ref{est1}}]
Consider $\Sp(k)$ with $k\ge 2$.
Denote by $\varepsilon_i- \varepsilon_j$,
$\pm(\varepsilon_i+\varepsilon_{j})$ for $1\le i,j\le k$, $i\ne
j$, and $\pm2\varepsilon_{i}$ for $1\le i\le k$ the roots of the
Lie algebra $\frak{sp}(k)$. Its positive roots
are $\varepsilon_i\pm\varepsilon_j$ for $1\le i<j\le k$ and
$2\varepsilon_i$.
The fundamental weights are
$$
\omega_1=\varepsilon_1,\ \omega_2=\varepsilon_1+\varepsilon_2,
\dots,\omega_{k}=\varepsilon_1+\varepsilon_2+\dots+
\varepsilon_{k}.
$$
Let $\omega=\sum_{i=1}^{k}m_i\omega_i$ be the dominant weight of
a complex irreducible representation which acts on a complex
vector space $V(\omega)$. The only complex irreducible
representations which are complexifications of
real representations are those for which
$$\sum_{i=0}^{\frac{k-1}{2}}m_{2i+1}$$
is even. (See \cite{S}, Theorem E, page 140.)
The other complex irreducible representations $V(\omega)$  determine
real irreducible representations of real dimension 
$2\dim_{\Cc} V(\omega)$.

$R(\Sp(k))$ is a polynomial ring in exterior powers which are self-conjugate.
A complex representation is real if and only if 
it is represented by a polynomial with coefficients as specified in the proof 
of Lemma \ref{ext-quaternionic}.

Put
$\omega+\delta=\sum_{i=1}^{k}g_i\varepsilon_i$.
In this case the Weyl Dimension Formula gives
$$\dim_{\Cc} V(\omega)=\prod_{1\le i<j\le k}\frac{(g_i-g_j)(g_i+g_j)}
{(j-i)(2k+2-i-j)}\prod_{1\le i\le k}\frac{g_i}{k-i+1}.$$
Since the inequality (\ref{eneq}) again holds, it is sufficient to compute
only the dimensions:
\begin{align*}
\dim_{\Cc} V(\omega_{j})&=\binom{2k+1}{j}\frac{2k-2j+2}{2k-j+2},
\quad 1\le j\le k,\\
\dim_{\Cc} V(2\omega_{1})&=(2k+1)k.
\end{align*}
For $k\ge 3$ consider a real irreducible representation $\rho$
which is determined by a complex irreducible representation with a dominant
weight $\omega$ different from $\omega_1$. Its dimension is at least
\begin{equation*}
\min\{\dim_{\Cc} V(\omega_2),
2\dim_{\Cc} V(\omega_3),\dim_{\Cc} V(2\omega_1)\}= 2k^2-k-1.
\end{equation*}
\end{proof}

\section{Appendix}
 
In this appendix we give a self-contained proof of the result of Davis and 
Mahowald (\cite{DaMa}) that the generator of $\pi_{8k}(\SO)=\Zz /2$,
for $k\ge 1$, lifts to $\pi_{8k}(\SO(6))$ and of a similar result 
that the generator of $\pi_{8k+4}(\Sp)=\Zz /2$ lifts to $\pi_{8k+4}(\Sp (1))$. 
These results are used in the proof of Lemma \ref{lift}.
At the very end we develop Remark \ref{note}
on possible generalizations of Theorem \ref{main}.

Let us write $C_n$ for the cofibre of the map
$z\mapsto z^n$: $S^1\to S^1$ (where $S^1$ is the space of complex
numbers of modulus $1$). Thus we have a cofibre sequence:
$$
S^1 \overset{n}\longrightarrow S^1 \longrightarrow C_n  \longrightarrow
S^2 \overset{n}\longrightarrow S^2 \to \cdots
$$

For pointed finite complexes $X$ and $Y$, we write
$\omega^0\{ X;\, Y\}$ for the group of stable maps $X\to Y$
and $\omega^i\{ X;\, Y\}$, with cohomological indexing,
for $\omega^0\{ X;\, S^i\wedge Y\}$.
In $K$-theory we write in the same way
$KO^0\{ X;\, Y\}$ for $[{\bf X};\, {\bf Y}\wedge {\bf KO}]$,
where ${\bf X}$ and ${\bf Y}$ are the suspension spectra of $X$ and $Y$, and 
${\bf KO}$ is the real $K$-theory spectrum,
and $KO^i\{ X;\, Y\} = KO^0\{ X;\, S^i\wedge Y\}$.
There is a Hurewicz map (or $d$-invariant)
$$
\omega^i\{ X;\, Y\} \to KO^i\{ X;\, Y\} .
$$
The groups $\omega^i\{ X;\, Y\}$ and $KO^i\{ X;\, Y\}$ 
can be identified with $\widetilde \omega^i(X\wedge D(Y))$,
and $\widetilde KO^i(X\wedge D(Y))$, respectively,
where $D(Y)$ is the stable dual of $Y$.

We shall also need the cohomology theory $J^*$ defined to be the
fibre of the stable Adams operation $\psi^3-1$ as a self-map of the
the real $KO$-theory spectrum localized at $2$. It is 
thus related to $KO$-theory by a long exact sequence:
$$
\cdots\to
KO^{*-1}_{(2)} \overset{\psi^3-1}\longrightarrow KO^{*-1}_{(2)} \to J^*
\to KO^*_{(2)} \overset{\psi^3-1}\longrightarrow  KO^*_{(2)} \to \cdots
$$

For any space $X$ and integer $k$,
the smash product with the identity on $X$ gives a map
$$
\wedge\, 1_X: \Zz \cong KO^{-8k}(*) \to KO^{-8k}\{ X;\, X\}\, .
$$
The image of one of the generators will be called
a {\it Bott element}.
Suppose that the rational homology of $X$ is zero.
We call a stable map $A\in \omega^{-8k}\{ X;\, X\}$ an
{\it Adams map} if its Hurewicz image in $KO^{-8k}\{ X;\, X\}$
is a Bott element. (See \cite{CK1} and \cite{CK2}.)

The proof of the result of Davis and Mahowald is based on the existence
of an unstable Adams map on $C_8$.
\begin{prop} 
There is an unstable Adams map
$$
A : \Sigma^{15}C_8 \to \Sigma^{7} C_8 .
$$
More precisely, the Hurewicz map
$$
[ \Sigma^{15}C_8 ;\, \Sigma^{7} C_8 ] \to KO^{-8}\{ C_8;\, C_8\}
=\Zz /8\oplus \Zz /2
$$
is surjective.
\end{prop}

\begin{proof} 
We start from the standard fact that $\pi_{14}(S^7)=(\Zz /120)\xi$,  
where $\xi$ stabilizes to $2\sigma\in \omega_7(*)=(\Zz /240)\sigma$. 
From the diagram of exact sequences of the defining cofibration 
of $C_8$:
$$
\xymatrix{
{}\ar[r]^-{8} & [S^{15};\, S^7] \ar[d] \ar[r] 
&[\Sigma^{13}C_8;\, S^7] \ar[d] \ar[r]
&\Zz /120= [S^{14};\, S^7]\ar[r]^-{8} \ar[d] &{}\\
{}\ar[r]^-{8} & \omega_8\ar[r] \ar[d] 
&\widetilde\omega^{-6}(C_8)\ar[d] \ar[r] 
& \Zz /240=\omega_7 \ar[r]^-{8} \ar[d]^-{h_J} & {}\\
{}\ar[r]^-{8} & J_8\ar[r] \ar[d] 
&\widetilde J^{-6}(C_8)\ar[d] \ar[r] 
& \Zz /16=J_7 \ar[r]^-{8} \ar[d] & {}\\
{}\ar[r]^-{8} & \Zz_{(2)}=(KO_8)_{(2)}\ar[r] 
&\Zz /8=\widetilde KO^{-6}(C_8)_{(2)} \ar[r]
& 0=(KO_7)_{(2)} \ar[r]^-{8} & {}
}$$
we can see that there is an element $a\in [\Sigma^{13}C_8;\, S^7]$
lifting an odd multiple of $\xi$ and mapping to the generator of 
$\Zz /8=\widetilde KO^{-6}(C_8)$.
The reason is that the map $\widetilde J^{-6}(C_8) \to \widetilde 
KO^{-6}(C_8)$ is onto
(since the Adams operator $\psi^3$ acts as multiplication by $3^4$ on
$(KO_8)_{(2)}$, and so $\psi^3-1$ as $80$ on $\widetilde KO^{-6}(C_8)_{(2)}$)
and that the Hurewicz homomorphism $h_J:\omega_7\to J_7$ is also onto,
which follows from the computation of the $J$-homomorphism and the 
$e$-invariant in \cite{Adams2}, Theorem 1.6,
and the relation between the $e$-invariant and $h_J$ (see \cite{CK3}, 
Section 1).

The Moore space $C_n$ is essentially self-dual. For we have a cofibration
sequence:
$$
D(S^1) =S^{-1} \overset{n}\longleftarrow D(S^1)=S^{-1}
\longleftarrow D(C_n)
\longleftarrow D(S^2)=S^{-2}  \overset{n}\longleftarrow D(S^2)=S^{-2},
$$
which allows us to identify the stable dual
$D(C_n)$ with $\Sigma^{-3}C_n$. 
The stable duality is specified by two structure maps
$$
S^3 \to C_n\wedge C_n \to S^3.
$$
We choose an unstable representative
$$
d: S^{15} \to C_8 \wedge S^{12}\wedge C_8
$$
of the first structure map.

We shall also need the standard stable equivalence:
$$
C_8\wedge C_8 \cong \Sigma^1C_8 \vee \Sigma^2C_8.
$$
This comes from the basic cofibration sequence by smashing with $C_8$:
$$
S^1 \wedge C_8 \overset{8}\longrightarrow S^1\wedge C_8 \longrightarrow 
C_8\wedge C_8
\longrightarrow S^2\wedge C_8 \overset{8}\longrightarrow S^2\wedge C_8,
$$
because $8=0\in\omega^0\{ C_8;\, C_8\}$.
We choose an unstable representative:
$$
p: \Sigma^{12} C_8 \wedge C_8 \to \Sigma^{13} C_8
$$
of the projection onto the first factor.

Now we can define 
the Adams map $A$ (modulo slight adjustment which will be described
later) as the composition:
$$
C_8 \wedge S^{15} \xrightarrow{1\wedge d}
C_8 \wedge S^{12}\wedge C_8 \wedge C_8 
\xrightarrow{p\wedge 1} S^{13}\wedge C_8 \wedge C_8
\xrightarrow{a\wedge 1}{}\, S^7\wedge C_8.
$$
We have to check that its Hurewicz image in $KO$-theory is a Bott element.
We have
\begin{align*}
KO^{-8}\{ C_8;\, C_8\}
&=\widetilde KO^{-5}(C_8\wedge C_8)
=\widetilde KO^{-6}(C_8)\oplus \widetilde KO^{-7}(C_8),\\
\omega^{-8}\{ C_8;\, C_8\}
&=\widetilde \omega^{-5}(C_8\wedge C_8)
=\widetilde \omega^{-6}(C_8)\oplus \widetilde \omega^{-7}(C_8)
\end{align*}
by duality and the stable splitting.
From the cofibration exact sequence:
$$
\xrightarrow{\,8\,} KO_{9}=\Zz /2 \to \widetilde KO^{-7}(C_8) \to 
KO_{8}=\Zz\xrightarrow{\,8\,} KO_8=\Zz \to \widetilde 
KO^{-6}(C_8) \to KO_{7}=0
$$
we see that
$KO^{-8}\{ C_8;\, C_8\} =\Zz /8\oplus\Zz /2$.
Indeed, the same calculation determines the stable homotopy:
$\omega^0\{ C_8;\, C_8\}=\Zz /8\oplus\Zz /2$ generated 
by the identity map of
order $8$ and an element of order $2$ given by the composition:
$$
h:\Sigma C_8 \to S^3 \xrightarrow{\,\eta\,} S^2 \to \Sigma C_8
$$
of the maps in the cofibration sequence and the Hopf element $\eta$.
The smash products of these maps with a Bott class $v\in KO_8$ give 
generators of $KO^{-8}\{ C_8;C_8\}$.

The calculation will be made by considering the commutative diagram
$$
\xymatrix{
[S^{13}\wedge C_8,\, S^7]\ar[r]\ar[d]
&\widetilde\omega^{-6}(C_8)\ar[r]^-{h_{KO}}\ar[d]^-{i_1}
&\widetilde KO^{-6}(C_8)\ar[d]^-{i_1}\\
[S^{15}\wedge C_8,\, S^7\wedge C_8]\ar[r]
&\omega^{-8}\{C_8;\, C_8\}\ar[r]^-{h_{KO}}
&KO^{-8}\{C_8;\, C_8\}
}
$$
in which 
the first vertical map is given by the smash with identity on $C_8$ 
composed with $(p\wedge 1)^*$ and $(1\wedge d)^*$, 
$i_1$ is the inclusion of the first summand and 
$h_{KO}$ is the Hurewicz homomorphism. 
The element $a\in [S^{13}\wedge C_8;\, S^7]$ maps to $A\in
[S^{15}\wedge C_8;\, S^7\wedge C_8]$ and simultaneously to
a generator of $\widetilde KO^{-6}(C_8)=\Zz /8$ as shown above. 
This implies that the Hurewicz image of $A$ is equal to
$v\cdot (m+nh )$, where $m$ is an odd integer, $n\in \Zz /2$,
and $v\in KO_8$ is a Bott class. We can certainly multiply $A$ by
an odd integer to arrange that $m=1$.
We can also modify $A$ by the element $\Sigma^6h$ if necessary to achieve
$n=0$, because $h^2=0$ (stably).
These adjustments produce the required unstable
Adams map on $C_8$. We have noted in passing that $vh$ lies in the
image of the Hurewicz homomorphism, and this completes the proof.
\end{proof}

Now consider the diagram:
$$\xymatrix{
{}\ar[r]^-{8} &\widetilde KO^{-1}(S^9)=\Zz /2 \ar[r]
&\widetilde KO^{-1}(\Sigma^7C_8)\ar[r] 
&\widetilde KO^{-1}(S^8)=\Zz /2 \ar[r]^-{8} &{}\\
{}\ar[r]^-{2} &\widetilde KO^{-1}(S^9)=\Zz /2 \ar[r]\ar[u]^-{4}_-{=0}  
&\widetilde KO^{-1}(\Sigma^7C_2)\ar[r]\ar[u]
&\widetilde KO^{-1}(S^8)=\Zz /2 \ar[r]^-{2} \ar[u]^-{1}&{}
}$$
The generator $y\in \pi_8(O)=\widetilde KO^{-1}(S^8)=\Zz /2$ has
a lift, $x$ say, a generator of $\widetilde KO^{-1}(\Sigma^7 C_8)$
coming from a generator of $\widetilde KO^{-1}(\Sigma^7C_2)\cong \Zz /4$.
This class $x$ gives a map $\Sigma^7 C_8 \to \O$.
We shall show that it lifts to an element $\tilde x\in
[\Sigma^7C_8; \SO(6)]$.

Recall from obstruction theory that, if $X$ is a pointed
finite complex with $\dim X <2n-1$,
the obstruction to lifting a class $x\in [X; \O ]$
to $[X; \O (n)]$ is precisely the image of $\theta (x)$ in 
$\omega^0\{ X; P^\infty_n\}$ (where $P^\infty_n = P(\Rr^\infty )/P(\Rr^n)$
and $\theta$ is the stable splitting used in Section 3).
Since $\dim \Sigma^7C_8  = 9 < 2\cdot 6 -1$, it suffices to show that
the obstruction vanishes in
$\omega^0\{ \Sigma^7 C_8; P^\infty_6\}$.
We look at the diagram in stable homotopy corresponding to the
$KO$-theory diagram above:
$$\xymatrix{
{}\ar[r]^-{8} &\widetilde \omega_{9}(P^{\infty}_6)=\Zz /12 \ar[r]  
& \omega^{0}\{ \Sigma^7C_8;\, P^{\infty}_6\}\ar[r]
&\widetilde \omega_{8}(P^{\infty}_6)=0 \ar[r]^-{8} &{}\\
{}\ar[r]^-{2} &\widetilde \omega_{9}(P^{\infty}_6)=\Zz /12 
\ar[r]\ar[u]^-{4}
&\omega^{0}\{ \Sigma^7C_2;\,P^{\infty}_6\}\ar[r]\ar[u]
&\widetilde \omega_{8}(P^{\infty}_6)=0 \ar[r]^-{2} \ar[u]^-{1}&{}
}$$
(The calculations of the
stable homotopy groups $\widetilde\omega_9(P^{\infty}_6)=\pi_9(\O/\O(6))$ 
and $\widetilde\omega_8(P^{\infty}_6)$ $=\pi_8(\O/\O(6))$ can be found in 
\cite{Pa}.) 
   From the diagram, the map
$\omega^0\{ \Sigma^7C_2;\, P^{\infty}_6\} \to
\omega^0\{ \Sigma^7C_8;\, P^{\infty}_6\}$ is zero, and
so the obstruction to lifting $x$ to $\SO(6)$ is zero.

We can now use iterates of the Adams map $A$ to lift the images of
$x$ under Bott periodicity to $\SO(6)$ as $(A^{k-1})^*(\tilde x)$:
$$
\xymatrix{ 
[\Sigma^7C_8,\, \SO(6)] \ar[r]^-{A^*}\ar[d]
&[\Sigma^{15}C_8,\, \SO(6)] \ar[r]^-{(\Sigma^8A)^*}\ar[d]
&{\cdots}\\
\widetilde KO^{-1}(\Sigma^7C_8) \ar[r]^-{v}_-{\cong}\ar[d]
&\widetilde KO^{-1}(\Sigma^{15}C_8) \ar[r]^-{v}_-{\cong}\ar[d] 
&{\cdots}\\
\widetilde KO^{-1}(S^8) \ar[r]^-{v}_-{\cong}
&\widetilde KO^{-1}(S^{16}) \ar[r]^-{v}_-{\cong} 
&\cdots
}
$$
(where $v$ is the Bott periodicity class).
Each vertical composition 
$$
[\Sigma^{8k-1}C_8,\, \SO(6)]
\to \widetilde KO^{-1}(\Sigma^{8k-1}C_8)
\to \widetilde KO^{-1}(S^{8k})=\Zz /2
$$
is surjective (since the first one is) and factors through
$$
[\Sigma^{8k-1}C_8, \, \SO(6)]
\to [S^{8k},\, \SO(6)]
\to \widetilde KO^{-1}(S^{8k})=\Zz /2 .
$$
We have thus established:

\begin{prop}[Davis, Mahowald]\label{so6}
For $k\ge 1$, 
the generator of $\pi_{8k}(O)=\widetilde KO^{-1}(S^{8k})=\Zz /2$ lifts
to an element $($with order a divisor of $8)$
in $\pi_{8k}(\SO(6))$.
\end{prop}

\begin{cor}
The generator of 
$\pi_{8k+1}(O)=\widetilde KO^{-1}(S^{8k+1})=\Zz /2$ lifts
to an element of order $2$
in $\pi_{8k+1}(\SO(6))$.
\end{cor}

\begin{proof} 
A lift is obtained by composing
with the Hopf element $S^{8k+1}\to S^{8k}$, which has order $2$.
\end{proof}

\begin{prop}\label{sp1} The generator of $\pi_{8k+4}(\Sp)=
\widetilde KSp^{-1}(S^{8k+4})=\Zz /2$ lifts to an element in
$\pi_{8k+4}(\Sp(1))$ for $k\ge 1$.
\end{prop}

\begin{proof} According to \cite{MT1} and \cite{MT2}, p. 261, the homomorphism
$$\pi_{12}(\Sp(1))=\Zz /2\oplus \Zz /2\to\pi_{12}(\Sp)=\Zz /2$$
induced by the standard inclusion is an epimorphism.

By the Bockstein sequence any element of $\pi_{12}\big(\Sp(1)\big)$
can be factored through $\Sigma^{11} C_8$. Hence the composition
in the first column of the following diagram is surjective:
$$
\xymatrix@C=3pc@R=2pc{
\left[\Sigma^{11} C_8,\, \Sp(1)\right] \ar [d] \ar[r]^-{(\Sigma^4 A)^*}&
\left[\Sigma^{19} C_8,\, \Sp(1)\right] \ar [d] \ar[r]^-{(\Sigma^{12} A)^*}&
\cdots\\
\widetilde KSp^{-1}(\Sigma^{11} C_8) \ar [d]
\ar[r]^-{v}_-{\cong}&
\widetilde KSp^{-1}(\Sigma^{19} C_8) \ar [d]
\ar [r]^-{v}_-{\cong} &
\cdots\\
\widetilde KSp^{-1}(S^{12})  \ar[r]^-{v}_-\cong& 
\widetilde KSp^{1}(S^{20}) \ar[r]^-{v}_-{\cong}&\cdots
}
$$
Since $A$ induces a Bott isomorphism in $KO$-theory, it induces
an isomorphism in $KSp$-theory as well. Consequently, all
vertical compositions
$$
\left[ \Sigma ^{8k+3} C_8,\, \Sp(1)\right]\to \left[\Sigma^{8k+3}
C_8,\, \Sp\right]\to\left[ S^{8k+4}, \Sp\right]
$$
are surjective and factor through
$$
\left[ \Sigma ^{8k+3} C_8,\, \Sp(1)\right]\to \left[S^{8k+4},\,
\Sp(1)\right]\to  \left[ S^{8k+4},\, \Sp\right]\,.
$$
\end{proof}

\begin{remark} 
Let $M$ be a connected, stably parallelizable, closed manifold of odd dimension
$n\ge 9$. There is, up to isomorphism, exactly one stably trivial, but not 
trivial, $n$-dimensional vector bundle $\tau$ over $M$. It may be described
as the pullback of the tangent bundle of $S^n$ by a map $f:M\to S^n$ of degree 
one (collapsing the complement of an open disc).  

In \cite{Th} Thomas showed that the structure group of $\tau$ can be reduced
to $SO(k)$ by the standard inclusion $SO(k)\hookrightarrow SO(n)$ if and only 
if $m\equiv 0$ mod~$a(m-k)$, where $n=m-1$ as in Theorem \ref{main}(A).
(The special case in which $\tau$ is the tangent bundle of $M$ was considered 
in \cite{BK}.) One might ask whether Theorem \ref{main} admits a 
generalization to such a bundle $\tau$ over $M$.

Consider a homomorphism  $\rho:G_k\to \O(n)$ for which there is a map
$\rho'':BG_{\infty}\to BO$ such that 
$\rho''\circ\iota_k\simeq\iota_n\circ B\rho$:
$$
\xymatrix{
BG_k \ar[r]^{B\rho}\ar[d]_{\iota_k}  & B\O (n) \ar[d]^{\iota_n}\\
BG_{\infty} \ar[r]^{\rho''}  & B\O
}
$$
and suppose that $\tau:M\to BO(n)$ can be factored through $B\rho$ by a map
$\xi:M\to BG_k$ such that $\tau=B\rho\circ\xi$. Now if, as in Lemma \ref{lift},
there is a map $\xi'$ such that $B\rho\circ\xi'=\tau$ and 
$\iota_k\circ\xi':M\to BG_{\infty}$ is trivial, then we may argue as in the 
proof of (A), (B) and (C) of Theorem \ref{main}.
Since $\xi'$ can be lifted to $\hat\xi:M\to G_{\infty}/G_k$
and a stable parallelization
of $M$ gives a stable splitting 
$h:S^n\to M$ of the map $f:M\to S^n$, we may take  
$x=\hat\zeta\circ h\in \widetilde\omega_n(G_{\infty}/G_k)$
satisfying the assumptions of Proposition \ref{hurewicz}. Consequently, 
$\rho:G_k\to G_{\infty}$ is described by (A), (B) and (C) of 
Theorem \ref{main}.

In the case that $\rho:\SO (k)\to \O (n)$ is the standard inclusion
we may take $\rho''$ to be identity, $\xi=\xi'$ satisfies the required 
condition and we reprove the result of \cite{Th}.

There are other special cases where the existence of $\xi'$ is guaranteed:
for example if $G_k=SU(k)$ and the map $K^o(M)\to KO^o(M)$ is injective,
as happens when $M$ is a product of spheres of the form 
$S^{4i_1}\times\dots\times S^{4i_l}\times S^{2j-1}$.
\end{remark}

\begin{ack}
Work on this paper started when the first author visited the Department
of Mathematical Sciences of the University of Aberdeen and was
completed when the second author visited the Department of Algebra and 
Geometry of the Masaryk University of Brno. The financial
support of the London Mathematical Society and the hospitality of 
both Departments are gratefully acknowledged.
The authors would also like to thank Mamoru Mimura for drawing their attention 
to this problem.
\end{ack}

\end{document}